\documentclass[10pt]{article}
\usepackage[latin1]{inputenc}
\usepackage[T1]{fontenc}
\usepackage{amsmath,amssymb,amstext}
\usepackage{theorem}
\usepackage{a4wide}
\usepackage{multicol}
\usepackage[english]{babel}
\usepackage{enumerate}
\newcommand{\R}{\mathbb{R}}
\newcommand{\N}{\mathbb{N}}

\newcommand{\E}{\mathbb{E}}

\newcommand{\dist}{\overset{(d)}{=}}
\newtheorem{thm}{Theorem}[section]

\newtheorem{lem}[thm]{Lemma}

\newtheorem{prop}[thm]{Proposition}

\theoremstyle{definition}
\newtheorem{defn}[thm]{Definition}
\newtheorem{exmp}[thm]{Example}
\newtheorem{rem}[thm]{Remark}

\newenvironment{dem}{\textbf{Proof}~:\\ }{\flushright$\blacksquare$\\}

\title{On processes which are infinitely divisible with respect to time}
\author{Roger \textsc{Mansuy}\thanks{\textbf{mansuy@ccr.jussieu.fr}, Laboratoire de Probabilit\'es et Mod\`eles Al\'eatoires,
Universit\'e Pierre et Marie Curie-Paris VI, Casier 188, 4 place
Jussieu, F-75252 Paris Cedex 05}}
\date{\today}
\begin{document}
\maketitle
\begin{abstract}The aim of this short note is to present the notion of IDT processes, which is a
wide generalization of Lévy processes obtained from a modified
infinitely divisible property. Special attention is put on a
number of examples, in order to clarify how much the IDT processes
either differ from, or resemble to, Lévy processes.
\end{abstract}
\noindent\textbf{Keywords~:} \emph{Infinitely divisible law, Lévy processes, Gaussian processes.}\\
\textbf{AMS 2000 subject classification ~:}*60G48, 60G51, 60G44,
60G10
\section{Introduction}
Motivated by some one-to-one map from the set of infinitely
divisible laws, or rather from the set of their Lévy measures onto
itself which was noticed by Barndorff-Nielsen and Thornbj\o rnsen
\cite{MR2003j:46101a}, it appears to be of some interest to
consider the class of stochastic processes $(X_t;t\geq 0)$ which
enjoy the following property~:
\begin{eqnarray}\forall n\in \N^*,\
(X_{nt}; t\geq 0)&\dist&(X^{(1)}_t+...+X^{(n)}_t;t\geq 0)
\label{idt}\end{eqnarray}where $X^{(1)}$,..., $X^{(n)}$ are
independent copies of $X$.\\ We shall call such a process an IDT
process (which stands for \textbf{I}nfinitely \textbf{D}ivisible
with respect to \textbf{T}ime). A simple remark, whose proof is
left to the reader, is
\begin{prop}:\\
Any Lévy process is IDT
\end{prop}

However, there are many other processes than Lévy processes which
are IDT, and the purpose of this paper is to exhibit large sets of
IDT processes. More precisely, in Section 2, we give some general
examples of IDT processes; in Section 3, we characterize the
Gaussian processes which are IDT. In Section 4, we focus on the
resemblance between a given IDT process and the Lévy process with
the same one-dimensional marginals at fixed times. Section 5 is
devoted to a description of IDT processes as particular
path-valued infinitely divisible variables.

\section{Some general examples of IDT processes}
The aim of this section is to provide as many examples as possible
of IDT processes  which are not Lévy processes.
\begin{paragraph}{(2.1)} Here is a first simple example; given a
strictly stable random variable $S_{\alpha}$ with parameter
$\alpha$, the process $X$ defined by $$
X_t=t^{1/\alpha}S_{\alpha},\ t\in\R^+,$$ is an IDT process. This
simply relies on the fact that $n^{1/\alpha}S_{\alpha}$ is equal
in law to the sum of $n$ independent copies of $S_{\alpha}$.
\end{paragraph}
\begin{paragraph}{(2.2)}A more interesting set of examples is
obtained by transforming a Lévy process, more generally an IDT
process, by integration in time with respect to a measure. Namely,
if $X$ is an IDT process and $\mu$ a measure on $\R^+$ such that
$$\ X^{(\mu)}_t=\int\mu(du)X_{ut},\ t\in\R^+,$$ is well
defined, then $X^{(\mu)}$
is again an IDT process: this property follows easily from (\ref{idt}).\\
Some particular cases of interest may be obtained with
$\mu(du)=\frac{du}{u}1_{[a,b]}(u)$ for any $0<a<b<\infty$. In such
a case, $X^{(\mu)}_t=\int_{at}^{bt}\frac{du}{u}X_u$; in fact, this
is the original example in \cite{MR2003j:46101a} , which motivated
us to study IDT processes.\\Thus, we note that this procedure
allows  to construct some quite regular IDT processes, i.e. if $X$
is a Lévy process,then $\int_0^t\frac{du}{u}X_u$ is an absolutely
continuous IDT process.
\end{paragraph}
\begin{paragraph}{(2.3)} Some similar examples can be provided using
only the jumps of an IDT process $X$.\\ Indeed, if
$f:(0,\infty)\times\R\rightarrow\R$ satisfies $f(.,0)\equiv0$,
then
$$X_t^{f}=\sum\limits_{v\geq 0}f(v,(\Delta X)_{vt})$$ is an IDT
process (this easily follows from the fact that, if $(X_t;t\geq
0)$ is an IDT process then the same holds for $((\Delta X)_t;t\geq
0)$), provided it exists.
\end{paragraph}
\section{IDT Gaussian processes}One may wonder which among centered Gaussian processes
$(G_t;t\geq 0)$ (which, for simplicity, we assume to be centered)
are IDT. In order to characterize IDT Gaussian processes, we first
recall the Lamperti transformation concerning self-similar
processes (see the original paper \cite{MR0138128} or
\cite{MR1920153} Theorem 1.5.1 p11)~:
\begin{lem}\label{lamperti}:\\ A process $(X_t;\ t\geq0)$ enjoys the scaling
property of order $h$, i.e. $$\forall \alpha>0,\qquad (X_{\alpha
t};\ t\geq0)\dist(\alpha^h X_{t};\ t\geq0),$$ if, and only if, its
Lamperti transform $(\tilde{X}_y:=e^{-hy}X_{e^y};\ y\in\R)$ is a
strictly stationary process.
\end{lem}
The next proposition emphasizes again how much more general IDT
processes are than Lévy processes (since the only Gaussian Lévy
processes are Brownian motions with drifts).
\begin{prop}\label{gaussidt}:\\ Let  $(G_t;\ t\geq 0)$ be a centered Gaussian process, which is assumed to be continuous in
probability (which is equivalent to its covariance function $c$
being continuous). Then the following properties are equivalent~:
\begin{enumerate}
    \item $(G_t;\ t\geq 0)$ is an IDT process.
    \item The covariance function $c(s,t):=\E[G_sG_t]$, $0\leq
    s\leq t$, satisfies $$\forall \alpha>0,\qquad c(\alpha
    s,\alpha t)= \alpha c(s,t),\qquad 0\leq s\leq t$$
    \item The process $(G_t;\ t\geq 0)$ satisfies the "Brownian
    scaling property", namely $$\forall \alpha>0,\qquad (G_{\alpha t};t\geq 0)\dist(\sqrt{\alpha}G_t;t\geq 0)$$
    \item The process $(\tilde{G}_y:=e^{-y/2}G_{e^y};\ y\in\R)$
    is stationary.
    \item The covariance fonction
    $\tilde{c}(y,z):=\E[\tilde{G}_y\tilde{G}_z]$, $y,z\in\R$, is
    of the form $$\tilde{c}(y,z)=\int \mu(du)e^{iu|y-z|},\qquad
    y,z\in\R$$where $\mu$ is a positive, finite, symmetric measure
    on $\R$.
\end{enumerate}
Then, under these equivalent conditions, the covariance function
$c$ of $(G_t;\ t\geq0)$ is given by
$$c(s,t)=\sqrt{st}\int\mu(da)e^{ia|\ln(\frac{s}{t})|}$$
\end{prop}

\begin{rem}(Private communication from F. Hirsch to M. Yor):\\The conditions found in Proposition \ref{gaussidt} are also
equivalent to the positivity of the quadratic form $Q$ defined as
$$Q(f)=\int_0^1du\,f(u)\int_0^1ds\,f(su) c(1,s),\qquad f\in
L^2([0,1])$$
\end{rem}
\begin{dem}
\begin{description}
    \item[$(1\Leftrightarrow2)$] The IDT property \eqref{idt} is
    equivalent to $$\forall n\in \N,\qquad c(ns,nt)=nc(s,t),\qquad
    0\leq s\leq t$$ and also to
    $$\forall q\in \mathbb{Q}_+,\qquad c(qs,qt)=qc(s,t),\qquad
    0\leq s\leq t$$
    The result is then obtained using the density
    of $\mathbb{Q}_+$ in $\R_+$ and the continuity of $c$ (deduced
    from the continuity in probability, hence in $L^2$, of $(G_t;\ t\geq0)$).
    \item[$(2\Leftrightarrow3)$] Simple, since the law of a centered Gaussian process is determined by its
covariance function.
    \item[$(3\Leftrightarrow4)$] Lamperti's transformation (Lemma
    \ref{lamperti}) of order $h=1/2$.
    \item[$(4\Leftrightarrow5)$] Bochner's theorem for definite
    positive functions.
\end{description}
\end{dem}

\begin{exmp}:\\ Let $\varphi\in L^2(\R^+,du)$; then the process
$G^{(\varphi)}$ defined by
\begin{eqnarray}
\forall t> 0,\
G^{(\varphi)}_t&=&\int_0^{\infty}\varphi\left(\frac{u}{t}\right)dB_u,\qquad
and\qquad G^{(\varphi)}_0=0\end{eqnarray} is an IDT process.\\
Indeed, it suffices to note that $G^{(\varphi)}$ is a well defined
Gaussian process with covariance
function$$c(s,t)=s\int_0^{\infty}dv\
\varphi\left(\frac{s}{t}v\right)\varphi(v)=\sqrt{st}\int\mu(dy)e^{iy|\ln(\frac{s}{t})|},\qquad\qquad
s,\ t>0$$with $\int \mu(dy)e^{i
yx}=e^{-|x|/2}\int_0^{\infty}dv\,\varphi(e^{-|x|}v)\varphi(v)$,
$x\in\R$.\\
In particular, we can compute $\mu$ in the following simple
cases~:$$\begin{array}{|ccc|c|}\hline&&&\\
   \ &\varphi(x)&\quad & \quad\mu(dy)\quad \\&&&\\\hline\hline
&&&\\
%   &1_{x\leq1}&&\cfrac{2}{\pi(4y^2+1)}dy\\&&&\\\hline
%   &&&\\
   &
   \begin{array}{c}
     x^{-\alpha}1_{x\geq 1} \\
     \\
     x^{-\alpha}1_{x\leq 1} \\
   \end{array}&
   \left.\begin{array}{c}
     \mbox{with } \alpha> 1/2 \\ \\
     \mbox{with } \alpha< 1/2 \\
   \end{array}\right\}&\cfrac{1}{2\pi(y^2+(1/2-\alpha)^2)}dy  \\&&&\\\hline
   &&&\\
   %&&&\cfrac{1}{2\pi(y^2+(1/2-\alpha)^2)}dy  \\&&&\\\hline
   %&&&\\
   &(1-x)^{-\alpha}1_{x\leq 1}&\mbox{with } \alpha< 1/2\ & \
   \sum\limits_{n=0}^{\infty}\cfrac{\Gamma(1-\alpha)}{\Gamma(\alpha)}\
   \cfrac{\Gamma(n+\alpha)}{\Gamma(n+2-\alpha)}\ \cfrac{n+1/2}{\pi(y^2+(n+1/2)^2)}\ dy\ \\
   &&&\\
   \hline
\end{array}$$
More generally, if
$\varphi(x)=\cfrac1{x^{\alpha}}\left(\sum\limits_{n=0}^{\infty}a_n
x^n\right)$ with some sequence of positive numbers $(a_n;\
n\geq0)$, then $$\mu(dy)=\sum\limits_{n=0}^{\infty}a_n\int dv\
\varphi(v)\ v^{n-\alpha}\
\cfrac{n-\alpha+1/2}{\pi(y^2+(n-\alpha+1/2)^2)}\ dy$$ Therefore,
the quantity $\int \mu(dy)e^{i y\sqrt{2\lambda}}$ may be
interpreted as the Laplace transform of a positive random variable
$A$ and $\mu$ is the law of $\beta_A$, with $(\beta_t;\ t\geq0)$ a
standard Brownian motion.
\end{exmp}

\section{IDT and Lévy processes with the same marginals} First we
note that the laws of finite dimensional marginals of an IDT
process $X$ are infinitely divisible. In particular, for any fixed
$t$, the law of $X_t$ is infinitely divisible.
\begin{prop}\label{mimi}:\\Let $(X_t;t\geq 0)$ be a right-continuous IDT process.
Denote by $(\tilde{X}_t;t\geq 0)$ the unique Lévy process such
that$$X_1\dist \tilde{X}_1$$ Then $(X_t;t\geq 0)$ and
$(\tilde{X}_t;t\geq 0)$ have the same one-dimensional marginals,\\
i.e. for any  fixed $t\geq 0$,\begin{eqnarray}X_t&\dist&
\tilde{X}_t\label{eq}\end{eqnarray}
\end{prop}
\begin{dem}
Since, for any $k\in\N$, $X_k\dist X^{(1)}_1+...+X^{(k)}_1$, it
follows that $X_k\dist \tilde{X}_k$.\\Identity (\ref{eq}) can then
be obtained for any rational time since the $n$-th power of the
characteristic functions of $X_{k/n}$ and of $\tilde{X}_{k/n}$ are
equal and non-vanishing (because the laws of these variables are
infinitely divisible).\\ We conclude by using the right-continuity
of paths of both $X$ and $\tilde{X}$.
\end{dem}
\begin{rem}:\\ Proposition \ref{mimi} just states that IDT processes
"mimick" Lévy processes in the sense of \cite{MR87k:60147} who
exhibits other examples of various processes with the same
1-dimensional marginals. Further studies in this direction can be
found in \cite{MR1914701} or \cite{MR1785391}.
\end{rem}

We now illustrate Proposition \ref{mimi} by computing the Lévy
measure of $X^{(\mu)}_1=\int_0^{\infty}d\mu(u)X_u$, where $X$ is a
Lévy process (hence $X^{(\mu)}$ is an IDT process; see (2.2)), or
equivalently of the Lévy process $(\tilde{X}^{(\mu)}_t, t\geq 0)$
related to $X^{(\mu)}$ by (\ref{eq}).
\begin{prop}:\\ Let $X$ be a Lévy process and $\nu$ its Lévy measure.
Let $\nu^{(\mu)}$ be the Lévy measure of the infinitely divisible
variable $X^{(\mu)}_1=\int_0^{\infty}\mu(du)X_u$ for a "good"
measure $\mu$. Then for any non-negative Borel function $f$, we
have\begin{eqnarray}
\int\nu^{(\mu)}(dy)f(y)&=&\int_0^{\infty}dh\int\nu(dx)f(\mu([h,\infty))x)\label{nu}\end{eqnarray}
\end{prop}
\begin{dem}First, integration by part implies
$$\int_0^{\infty}\mu(du)X_u=\int_0^{\infty}\mu([h,\infty))dX_h$$
so that
\begin{eqnarray*}\E\left[\exp\left(-\lambda\int_0^{\infty}\mu(du)X_u\right)\right]
&=&\exp\left(-\int_0^{\infty}dh\int\nu(dx)\left(1-e^{-\lambda\mu([h,\infty))x}\right)\right)
\end{eqnarray*} from which we immediately deduce (\ref{nu}).
\end{dem}
\begin{exmp}:\\With $\mu(du)=\frac{du}{u}1_{[0,1]}(u)$, the identity
(\ref{nu}) becomes~:
$$\nu^{(\mu)}(dv)=\left(\int\frac{\nu(dx)}{x}e^{-v/x}\right).dv$$
In particular, for $\nu(dx)=\frac{dx}{\sqrt{2\pi x}}1_{x>0}$, then
$\nu^{(\mu)}(dv)=\frac{dv}{\sqrt{2v}}1_{v>0}$.
\end{exmp}

\section{A link with path-valued Lévy processes} It follows from
the very definition of an IDT process that, viewed as a random
variable taking values in path-space (i.e. $D=D([0,\infty))$, the
space of right continuous paths over $[0,\infty)$), this random
variable, which we shall denote by $\bar{X}$ is infinitely
divisible. The Lévy-Khintchine representation theorem for such
variables has been discussed in
\cite{P-Y2}, \cite{MR58:2949} or \cite{MR87m:60018} Chapter 5, among others.\\
%Before investigating what do such Lévy-Khintchine representations
%look like,
First of all, we remark that there exist $D$-valued infinitely
divisible variables which are not IDT. For example, consider the
random variable $R$ associated with the path of a squared Bessel
process of dimension 1. Although the distribution of $R$ is known
to be infinitely divisible (See \cite{MR51:4433} or
\cite{MR2000h:60050} theorem 1.2 p440), it is not the law of an
IDT process. Indeed, if it were so, the identity~(\ref{idt})
combined with the scaling property of Bessel processes would
entail that a squared Bessel process of dimension $n$ is the
square of a one-dimensional Brownian motion multiplied by $n$, which is absurd.\\
Now, we try to understand some properties of IDT Lévy measures. To
avoid some confusion, the expression "IDT Lévy measure" will
indicate the Lévy measure of the IDT process considered as an
infinitely divisible $D$-valued variable (not the Lévy measure of
the mimicked Lévy process as in Proposition \ref{mimi}).
\begin{lem}\label{lemme}:\\ Let $X$ be a $D$-valued infinitely divisible
variable. $X$ is an IDT process if, and only if,
\begin{itemize}
    \item its Lévy measure $M$ over $D$ satisfies, for any
non-negative functional $F$ on $D$,\begin{eqnarray} \int_D
M(dy)F(y(n \cdot))
&=&n\int_DM(dy)F(y(\cdot))\label{machin}\end{eqnarray}
    \item its Gaussian measure $\rho$ over $D$ satisfies, for any
    $t,u\in\R$,\begin{eqnarray}
\int_Dy(nu)y(nt)\rho(dy)&=&n\int_Dy(u)y(t)\rho(dy)
    \label{machin2}\end{eqnarray}
\end{itemize}
\end{lem}
\begin{dem}The IDT property (\ref{idt}) admits the following equivalent formulation~: for
any $f\in C_c(\R_+, \R_+)$,
\begin{eqnarray*}
\E\left[\exp{-\int_0^{\infty}dt f(t)X_{nt}}
\right]&=&\left(\E\left[\exp{-\int_0^{\infty}dt
f(t)X_{t}}\right]\right)^n
\end{eqnarray*}
That is
\begin{eqnarray*}
\int_D M(dy)\left(1-e^{-\int_0^{\infty}dt f(t)y(nt)}
\right)&=&n\int_DM(dy)\left(1-e^{-\int_0^{\infty}dt
f(t)y(t)}\right)
\end{eqnarray*}
and \begin{eqnarray*}\int_D\int_0^{\infty}du
f(u)y(nu)\int_0^{\infty}dtf(t)y(nt)\rho(dy)
&=&n\int_D\int_0^{\infty}du
f(u)y(u)\int_0^{\infty}dtf(t)y(t)\rho(dy)
\end{eqnarray*} from which (\ref{machin})-\eqref{machin2} follow
in a straightforward manner.
\end{dem}
This lemma provides us with some "new" constructions of IDT
processes~:
\begin{prop}\label{prop}:\\
Let $N$ be a Lévy measure on path space $D$. Define $M$ as
follows\begin{eqnarray}
\int_DM(dy)F(y(\cdot))&=&\int_0^{\infty}du\int_DN(dy)F(y(\frac{\cdot}{u}))\end{eqnarray}
Then $M$ is an IDT Lévy measure.
\end{prop}
\begin{dem}From Lemma \ref{lemme}, all we need to show is that, for any $n\in\N$,
\begin{eqnarray*}n\int_D M(dy)F(y(\cdot))&=&\int_D M(dy)F(y(n\cdot))
\end{eqnarray*}
This follows from the obvious change of variable $nv=u$
\begin{eqnarray*}
\int_0^{\infty}du\int_D
N(dy)F(y(\frac{n}{u}\cdot))&=&n\int_0^{\infty}dv\int_D
N(dy)F(y(\frac{\cdot}{v}))
\end{eqnarray*}
\end{dem}
\begin{exmp}:\\ In order to illustrate this Proposition \ref{prop}, let us
compute the IDT Lévy measure of some IDT processes construted in
paragraph (2.2). Let $X$ be a subordinator without drift, $\nu$
its Lévy measure, $\varphi$ a regular function and define
\begin{eqnarray*} X^{(\varphi)}_t&=&\int_0^{\infty}du\
\varphi(u)X_{ut}\end{eqnarray*} Then its IDT Lévy measure
satisfies for any functional $F$ over $D$
\begin{eqnarray}
\int_DM^{(\varphi)}(dy)F(y(\cdot))&=&\int_0^{\infty}du\int
\nu(dx)F(x\Phi(\frac{u}{\cdot}))
\end{eqnarray}
where $\Phi$ is the tail of the integral of $\varphi$~:
$\Phi(u)=\int_u^{\infty}dv\varphi(v)$.\\
To make a close link with Proposition \ref{prop}, the measure $N$
is now the image of $\nu$ by $x\mapsto x\Phi\left(\frac1x\right)$.
\end{exmp}
%\section{Some further generalizations}
%It may be of interest to study the following reinforcement of the
%IDT property.
%\begin{defn}~\\ A process $(X_t;\ t\geq0)$ will be said to be
%strongly IDT if for any finite sequence of positive numbers
%$(c_1,..., c_n)$ such that $c_1+...+c_n=1$, there exist some
%independent processes $(X^{(1)}_t;\ t\geq0)$... $(X^{(n)}_t;\
%t\geq0)$ such that $$ (X_t;\ t\geq0)\dist
%\left(X^{(1)}_{c_1t}+...+X^{(n)}_{c_n t};\ t\geq0\right)$$
%\begin{prop}~\\
%\end{prop}
%\begin{prop}~\\
%Any right-continuous IDT process $(X_t;\ t\geq0)$ is temporally
%self-decomposable, i.e. for any $c\in(0,1)$, there exists a
%process $(Y_t;\ t\geq0)$ such that $$(X_t;\ t\geq0)\dist (X_{ct};\
%t\geq0)+(Y_t;\ t\geq0)$$ with evident independence assumption.
%\end{prop}
%begin{dem}
%For any $\frac{p}{q}\in\mathbb{Q}$, we deduce from the IDT
%property \eqref{idt}~:
%\begin{eqnarray*}
%X_{pt};\ t\geq0)&\dist&(X^{(1)}_{t};\
%\geq0)+(X^{(2)}_{t}+...+X^{(p)}_{t};\ t\geq0)\\
%X^{(1)}_{t};\ t\geq0)&\dist& (\tilde{X}^{(1)}_{t/q};\
%t\geq0)+(\tilde{X}^{(2)}_{t/q}+...+\tilde{X}^{(q)}_{t/q};\ t\geq0)
%\end{eqnarray*}
%Therefore
%\begin{eqnarray*}
%(X^{(1)}_t;\ t\geq0)&\dist&(X_{pt/q};\ t\geq0)+(Y_t;\ t\geq0)
%\end{eqnarray*}
%with $(Y_t;\
%t\geq0)=(\tilde{X}^{(2)}_{t/q}+...+\tilde{X}^{(q)}_{t/q}-X^{(2)}_{t/q}-...-X^{(p)}_{t/q};\
%t\geq0)$
%\end{dem}
\section{Links with temporal self-decomposability}
In \cite{satobnm}, Barndorff-Nielsen, Maejima and Sato
 introduce the notion of temporally self
decomposable processes. These processes turn out to be deeply
linked with IDT processes.\\
We first recall the definition of temporal self-decomposability as
presented in \cite{satobnm}~:
\begin{defn}~:\begin{itemize}\item A real-valued process $(X_t;\ t\geq0)$ is
temporally self decomposable of order 1 if for any $c\in(0,1)$,
there exists a process $(U^{(c)}_t;\ t\geq0)$, called the
$c$-residual (of $(X_t; \ t\geq0)$), such that $X$ and
$(X_{ct}+U^{(c)}_t;\ t\geq0)$, with the obvious independence
assumption, have the same finite dimensional marginals, i.e. are
identical in law.\item $(X_t;\ t\geq0)$ is temporally self
decomposable of order $n>1$ if for any $c\in(0,1)$, the
$c$-residual $(U^{(c)}_t;\ t\geq0)$ is temporally self
decomposable of order $n-1$.\item $(X_t;\ t\geq0)$ is temporally
self decomposable of infinite order if $(X_t;\ t\geq0)$ is
temporally self decomposable of order $n$ for every $n\in
\N^*$.\end{itemize}
\end{defn}
\begin{prop}~:\\
A right-continuous IDT process is temporally self decomposable of
infinite order.
\end{prop}
\begin{dem}
Let $(t_1,....,t_n)\in \R^n_+$ and $z=(z_1,...,z_n)\in\R^n$.\\
Consider  the characteristic function of the
$(t_1,....,t_n)$-marginal of a right-continuous IDT process
$(X_t;\ t\geq0)$
$$\hat{\mu}_{t_1,....,t_n}(z):=\E\left[\exp{\left(i\sum\limits_{j=1}^nz_jX_{t_j}\right)}\right]$$
For any $r\in\mathbb{Q}$, the IDT property implies
\begin{eqnarray}
\hat{\mu}_{t_1,....,t_n}(z)&=&\left(\hat{\mu}_{rt_1,....,rt_n}(z)\right)^{\frac1r}\label{idtsato}
\end{eqnarray}
In particular, we deduce that, for $c\in\mathbb{Q}\bigcap(0,1)$~:
\begin{eqnarray*}
\hat{\mu}_{t_1,....,t_n}(z)&=&\hat{\mu}_{ct_1,....,ct_n}(z)\left(\hat{\mu}_{ct_1,....,ct_n}(z)\right)^{\frac1c-1}\\
&=&\hat{\mu}_{ct_1,....,ct_n}(z)\hat{\mu}_{\frac{c}{c'}t_1,....,\frac{c}{c'}t_n}(z)\qquad\mbox{with
}\frac1{c'}=\frac1c-1
\end{eqnarray*}
where the last equality is obtained using \eqref{idtsato}.\\
With the continuity assumption on $(X_t;\ t\geq0)$, we deduce that
$(X_t;\ t\geq0)$ is temporally self decomposable and that, for any
$c\in(0,1)$, the finite dimensional marginals of the associated
$c$-residual fit with the marginals of $X$ suitably rescaled;
hence the result.
\end{dem}
\section{Conclusion}
In this short note, we have introduced and discussed the notion of
IDT processes, and related this notion with the temporally self
decomposable processes (Section 6). We have shown that~:
$$\begin{array}{|c|}
  \hline\qquad\mathcal{L}\subset\mathcal{I}\subset\mathcal{S}_{\infty} \qquad\\\hline
\end{array} $$
Where $\mathcal{L}$ is the family of Lévy processes, $\mathcal{I}$
the family of IDT processes assumed to be right-continuous and
$\mathcal{S}_{\infty}$ the family of temporally self decomposable
processes of infinite order.

 \nocite{*}
\bibliographystyle{alpha}
\bibliography{idt}
\end{document}